\documentclass[12pt]{amsart}
\usepackage{amsfonts}
\usepackage{amsmath}
\usepackage{amssymb}
\usepackage[margin=1.2in]{geometry}
\usepackage{verbatim}
\newtheorem{theorem}{Theorem}[section]

\theoremstyle{remark}
\newtheorem{remark}[theorem]{Remark}

\newenvironment{acknowledgement}{%
	\par\medskip
	\noindent\textbf{Acknowledgements}%
}{\par}
\usepackage{tikz}
\usepackage{comment}
\usepackage[shortlabels]{enumitem}
\usepackage{graphicx}
\usepackage{float}
\usepackage[figurename=Fig.,labelfont=bf,labelsep=period]{caption}
\usepackage{array}
\newcolumntype{C}[1]{>{\centering\arraybackslash}m{#1}}
\usepackage{caption} 
\usepackage{tabularx,ragged2e,booktabs}
\newcolumntype{L}{>{\RaggedRight\arraybackslash}X}
\usepackage{adjustbox}
\makeatletter
\@namedef{subjclassname@2020}{%
  \textup{2020} Mathematics Subject Classification}
\makeatother
\usepackage{xcolor}
\usepackage[hidelinks]{hyperref}
\hypersetup{
    colorlinks = true,
    linkcolor = blue,
    anchorcolor = red,
    citecolor = blue,
    filecolor = red,
    urlcolor = red,
    linktoc=page
    }
\usepackage{cleveref}
\crefname{assumption}{Assumption}{Assumptions}
\usepackage{units}
\usepackage{setspace}
\newcommand{\ds}{\,\mathrm{d}s}

\newcommand{\abs}[1]{\left\lvert #1 \right\rvert}
\DeclareMathOperator{\Cont}{C}
\newcommand {\NN } {{\mathbb N}}
\usepackage{lineno} %for review
\begin{document}
\title[]{A Closed-form Approximation for Impulse Response of Fractionally Damped Oscillators}
\author[S.~Pathak]{Shashank Pathak$^1$}
\address[1]{School of Civil \& Environmental Engineering\\ Indian Institute of Technology (IIT) Mandi \\ Himachal Pradesh \\ India}
\author[M.~Ruzhansky]{Michael Ruzhansky$^2$}
\author[K.~Van~Bockstal]{Karel Van Bockstal$^2$} 
\address[2]{Ghent Analysis \& PDE center, Department of Mathematics: Analysis, Logic and Discrete Mathematics \\ Ghent University \\ 9000 Ghent \\ Belgium}
\email{shashank@iitmandi.ac.in, michael.ruzhansky@ugent.be, karel.vanbockstal@ugent.be}
\thanks{K. Van Bockstal is supported by the FWO Senior Research Grant G083525N}
\begin{abstract}  
We consider a fractionally damped oscillator, where the damping term is expressed by the Caputo fractional derivative of order $\beta\in (0,1).$ The impulse response of this oscillator can be expressed in terms of the bivariate Mittag-Leffler function consisting of a double infinite series. Although this series is uniformly convergent, its numerical implementation suffers from computational instabilities. In this contribution, we propose an approximate closed-form solution that avoids these numerical pitfalls while maintaining a reasonable accuracy. The resulting approximation is computationally efficient and robust, making it suitable for practical engineering applications. 
\end{abstract}
\maketitle
\tableofcontents
%%%%%%%%%%%
\section{Introduction}
%%%%%%%%
Modern structural engineering components use polymers, elastomers, viscoelastic dampers, and viscoelastic damping treatment for enhanced vibration resilience. These materials exhibit frequency-dependent and memory-effect (or hereditary) behavior. By now, it is well proven that the fractional (or non-integer) derivative-based models satisfactorily capture the mechanical behavior of systems involving such materials (see \cite{koh1990application, failla2020advanced, cunha2021new} and references therein). The fractional derivative models are compact, analytic, causal, and better fit the experimental data with a smaller number of calibration parameters compared to classical viscoelastic models and complex modulus methods \cite{bagley1979applications, gaul1989impulse}. 

The use of fractional material models leads to the fractional differential equation (FDE) of motion \cite{Pathak2024}. Comprehensive reviews detailing the analysis of these FDEs are available in \cite{rossikhin1997applications, rossikhin2010application}. A range of techniques is available to solve fractional systems, for example, the eigenvector expansion approach \cite{suarez1997eigenvector}, Galerkin approach \cite{das2013simple}, numerical methods \cite{padovan1987computational, singh2011algorithms}, and the integral transform (Laplace and Fourier transforms) based approaches \cite{kilbas2006theory, devillanova2016free, duan2021comparison}. Among these methods, the Laplace transform-based approaches are widely used by the engineering communities. This approach requires the inversion of the fractional-order transfer functions. Mostly, the studies employ Cauchy's residue theorem when evaluating the inverse transforms, which requires the numerical computation of complex integrals. Some methods, such as Oustaloup's method and the continued fraction approach, approximate fractional-order transfer functions by higher-order polynomials (see \cite{atherton2015methods} and references therein). However, higher-order polynomials lead to the loss of the intuitive understanding and computational simplicity.

Alternatively, the inverse transforms can also be expressed as an infinite series (involving Mittag-Leffler type functions) when differential equations have constant coefficients. However, such series solutions are generally not favoured in practical applications due to slow convergence, potential numerical instability (e.g., blow-up at large time-instants), and implementation difficulties
\cite{rossikhin2010application, duan2021comparison}. For example, Duan et al. \cite{duan2021comparison} attempted to investigate the impulse response of fractionally damped oscillators using the series solution. However, due to the slower convergence at higher time instants, they reverted to using the Laplace inversion method with complex integration via the Gauss–Laguerre quadrature scheme. 
  
\begin{figure}
\centering
\includegraphics[scale=0.4]{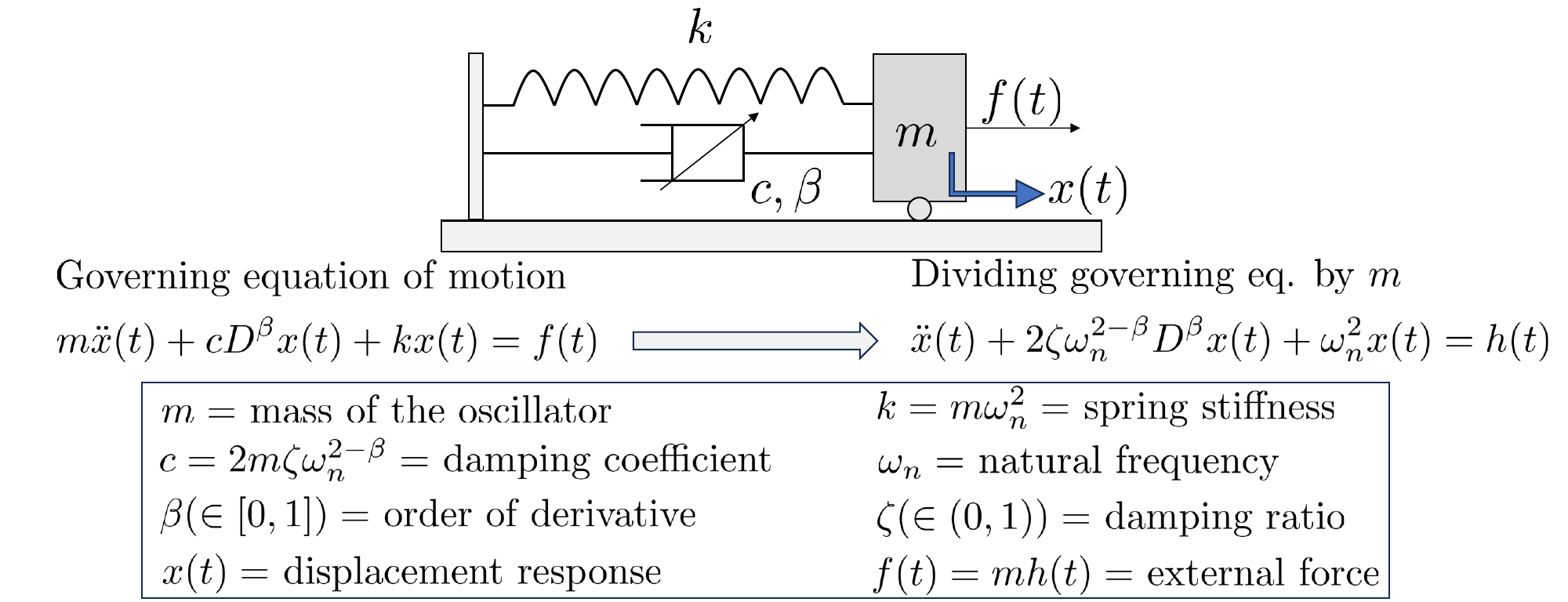}
\caption{A fractional single-degree-of-freedom oscillator along with the details of the governing equation of motion with zero initial conditions.}
\label{Fig:fsdf}
\end{figure}

In this contribution, we consider a fractionally damped oscillator as shown in \Cref{Fig:fsdf}. The fractional derivative is expressed by the Caputo fractional derivative $D^\beta$ of order $\beta \in (0,1)$ defined as \cite{caputo1967linear}:
$$
\left(D^\beta x\right)(t):=\left({}^{C}D_{0+}^{\beta}x\right)(t)=\frac{1}{\Gamma(1-\beta)}\int_0^t{(t-s)^{-\beta}x^\prime(s) \ds},$$
where $\Gamma(\cdot)$ is the gamma function and $x^\prime(s)$ is the first derivative of $x(s)$ with respect to $s$.

In \cite{Pathak2024}, a detailed derivation is presented that shows that the fractional material models lead to the fractional governing equation of motion 
\begin{equation}\label{eq:FDE}
\ddot{x}(t) + 2\zeta\omega_n^{2-\beta}  D^\beta x(t) + \omega_n^2	x(t) = h(t), \quad t >0,
\end{equation}
where $x(t)$ is the displacement response (unit: m) of the oscillator as a function of time $t$ (unit: s), $\ddot{x}(t)$ is the double derivative of $x(t)$ with respect to time, $\zeta\in(0,1)$ is the nondimensional damping ratio, and $\omega_n$ is the natural frequency of the oscillator (in rad/s). The function $h(t)$ is the applied excitation (acceleration in m/s$^2$) to the oscillator. The oscillator is initially at rest as shown in \Cref{Fig:fsdf}, i.e., the initial conditions are given by 
\begin{equation} \label{eq:ics}
x(0) = 0, \quad \dot{x}(0) = 0.
\end{equation}
It may be noted that $\omega_n$ has been raised to an exponent equal to $2-\beta$ due to dimensional consistency (as $D^
\beta x$ has unit m/s$^{\beta}$), see \cite{Pathak2024}. If $\beta=1$, this exponent becomes $1$, which is compatible with the classical viscously damped case, noting that $D^1 x(t) = \dot{x}(t)$. 

In this study, motivated by the aforementioned limitations of the series solutions, we will derive an approximate closed-form expression for the impulse response of fractionally damped systems governed by Equations \ref{eq:FDE} and \ref{eq:ics}. The proposed solution extends the classical impulse response of integer-order systems to the fractional case and can be convolved with arbitrary excitations to approximate the time-domain response of the system, making it suitable for practical engineering applications.

The paper is structured as follows. In \Cref{eq:series_solution}, we first derive the series solution to problem (\ref{eq:FDE}-\ref{eq:ics}). We then examine the limiting cases $\beta=0$ and $\beta=1$ in \Cref{subsec:limiting_case_beta0,subsec:limiting_case_beta1}. We show that the analytical summation of the resulting series coincides with the classical solution for the two extreme cases. 
Afterwards, in \Cref{subsec:blow_up_issues}, we discuss the numerical instabilities associated with the series solution when $0<\beta<1$ and propose a stable closed-form approximation in \Cref{subsec:approximate_solution}. The equivalent frequency and damping of this proposed solution are thoroughly discussed in \Cref{sec:equivalent_frequency,sec:equivalent_damping}. Moreover, in \Cref{sec:error_analysis}, we perform an error analysis between the proposed approximation and the series solution (or exact solution). Finally, in \Cref{sec:response}, we demonstrate the applicability of the developed closed-form equation for solving a fractionally damped system with arbitrary excitation. 
%%%%%%%%%%%%%
\section{Series Solution}
\label{eq:series_solution}
%%%%%%%%%%%%
To solve the problem (\ref{eq:FDE}-\ref{eq:ics}), we will first derive the impulse response of the oscillator and then, using the convolution integral, we shall obtain the desired response for the applied excitation $h(t)$.

The impulse response $I_\beta(t)$ is the solution to \eqref{eq:FDE} with $h(t)=0$ satisfying $x(0)=0$ and $\dot{x}(0) = 1. $ This solution can be obtained by applying the Laplace transform method. For more details, we refer the reader to \Cref{appendix:derivation_impulse_response}. It can be expressed in terms of the bivariate Mittag-Leffler function, which is defined by \cite[Eq. (3.16)]{Luchko1996}
\begin{equation}\label{eq:mlf}
E_{(a_1,a_2),b}(z_1,z_2) = \sum_{k=0}^\infty \sum_{\substack{ k_1+ k_2 = k\\ k_j\geq 0}} \binom{k}{k_1,k_2} \frac{\prod_{j=1}^2 z_j^{k_j}}{ \Gamma(b + \sum_{j=1}^2 a_j k_j)} =  \sum_{m=0}^{\infty} \sum_{l=0}^{\infty} \frac{\binom{m+l}{l}z_1^{l}z_2^{m}}{\Gamma(b + a_1 l + a_2 m)}.
\end{equation}
Here, \(\binom{k}{k_1, k_2} = \frac{k!}{k_1! k_2!}\) is the multinomial coefficient and  $\binom{m+l}{l}=\frac{(m+l)!}{m!l!}$ is the binomial coefficient. In particular, it holds that the impulse response $I_\beta(t)$ is given by 
\begin{align}\label{eq:imp_resp}
I_\beta(t) &= \sum_{m=0}^\infty \sum_{l=0}^\infty (-1)^{l+m}   \binom{l+m}{m}\left(2\zeta\omega_n^{2-\beta} \right)^m \omega_n^{2l} \frac{t^{2 l+(2-\beta)m + 1}}{ \Gamma(2l+ (2-\beta)m + 2) }  \\
&= t  E_{(2,2-\beta),2}(-(\omega_n t)^2,-2\zeta(\omega_n t)^{2-\beta}). \nonumber
\end{align}
Therefore, for any arbitrary fixed $T>0$, the solution $x$ to problem (\ref{eq:FDE}-\ref{eq:ics}) for given $h\in \Cont([0,T])$ belongs to $ \Cont^2([0,T])$ and is given by
%We will use a direct result provided in Corollary 4.9 of \cite{restrepo2021explicit} that is as follows:
%
\begin{equation}\label{eq:sol}
x(t) = \int_0^t h(t-s) I_\beta(s) \ds, \quad t>0, 
\end{equation}
see, e.g., \cite[Corollary~4.9]{restrepo2021explicit}. 
%%%%%%%%%%%%%%%%%%%%%%%%%%%%%%%%%%%%%%%%%%%%%%%%%%%%%%%%
\subsection{\texorpdfstring{$I_0(t)$}{}: Limiting case of \texorpdfstring{$I_\beta(t)$}{} when \texorpdfstring{$\beta=0$}{}}
\label{subsec:limiting_case_beta0}
%%%%%%%
For $\beta=0$, the governing equation becomes
$$
\ddot{x}(t) + (1+2\zeta)\omega_n^{2} {x}(t) = h(t).
$$
This equation does not have an explicit damping term. However, the damping ratio $\zeta$ contributes to the stiffness of the oscillator, causing the effective natural frequency $\omega_d$ to become 
\[\omega_d = \omega_n\sqrt{1+2\zeta}.\] From classical structural dynamics, we know that the impulse response of such an undamped ordinary oscillator is
\begin{equation}\label{eq:impulse_response_beta_0}
I_0(t) %=\dfrac{\sin{(\omega_n \sqrt{1+2\zeta})t}}{\omega_n\sqrt{1+2\zeta}} 
= \frac{\sin(\omega_d t)}{\omega_d}. 
\end{equation}
Using the series solution \Cref{eq:imp_resp} for $\beta =0,$ we have that
\[
I_0(t) = \sum_{k=0}^\infty (-1)^k \frac{t^{2k + 1 }}{ (2k + 1 )! } \sum_{m=0}^k  \binom{k}{m}\left(2\zeta \omega_n^2 \right)^m (\omega_n^2)^{k-m}, 
\]
which exactly matches the classical dynamics solution  \eqref{eq:impulse_response_beta_0}. For more details, the reader may refer to \Cref{appendix:beta_0}. 
%%%%%%%
\subsection{\texorpdfstring{$I_1(t)$}{}: Limiting case of \texorpdfstring{$I_\beta(t)$}{} when \texorpdfstring{$\beta=1$}{}}
\label{subsec:limiting_case_beta1}
%%%%%%%

For $\beta=1$, the governing equation takes the classical form (viscous damping)
$$
\ddot{x}(t) + 2\zeta\omega_n \dot{x} (t) + \omega_n^2 x (t) = h (t),
$$
whose impulse response is well-known and given by
\begin{equation}\label{eq:impulse_response_beta_1}
I_1(t)=e^{-\zeta\omega_n t}\dfrac{\sin{(\omega_d t)}}{\omega_d },    
\end{equation}
where \[\omega_d=\omega_n\sqrt{1-\zeta^2}\] denotes the effective or damped frequency (valid for $0<\zeta<1$).
Substituting $\beta=1$ into the series  \eqref{eq:imp_resp} yields
\begin{equation}\label{eq:limit_case:eq1}
I_1(t) = \sum_{m=0}^\infty \sum_{l=0}^\infty\frac{(l+m)!}{m!l!}\frac{(-1)^{l+m}\left(2\zeta\right)^m(\omega_nt)^{2l+m}t}{(2l+m+1)!},
\end{equation}
which recovers \eqref{eq:impulse_response_beta_1}. A rigorous derivation of this limiting case is provided in \Cref{appendix:beta_1}. 
%%%%%%%%%%%%%%%%%%%%%%%%%%%%%%%%%%%%%%%%%%%%%%%%%%%%%%%%
\subsection{Convergence and blow-up issues with \texorpdfstring{$I_\beta(t)$}{} when \texorpdfstring{$0<\beta<1$}{}}
\label{subsec:blow_up_issues}
\begin{figure}[h]
\centering
\includegraphics[scale=0.45]{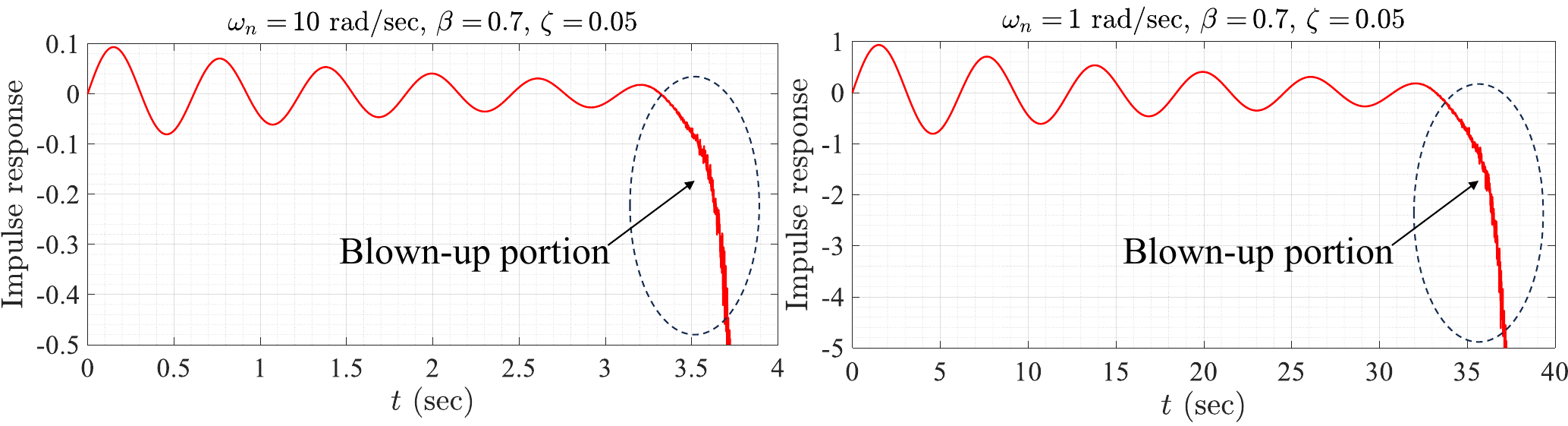}
\caption{Convergence issue in the series solution.}
\label{Fig:conv}
\end{figure}

%%%%%%%%%
Although the series \eqref{eq:imp_resp} is theoretically uniformly convergent (for $t \in [0,T]$ for arbitrary but fixed $T$), its numerical implementation suffers from significant computational instabilities. As illustrated in \Cref{Fig:conv}, if we take $\omega_n=10$ rad/s (let $\beta=0.7$ and $\zeta=0.05$), then we observe that the solution diverges (blows up) after approximately 3.5 seconds. This pattern scales with $\omega_n$, occurring after $t=35$ s if $\omega_n=1$ rad/s. 

For $\omega_n=10$ rad/s, at $t=3.5$ s, we need to sum at least till $m=13$ and $l=51$, in \cref{eq:imp_resp}, to obtain an acceptable accuracy. For these values of $m$ and $l$, we have in the denominator of \cref{eq:imp_resp}, $\Gamma(120.9)$ which is of the order of $4.143\times 10^{198}$ and in the numerator, we get $64!$ which is around $1.2689\times10^{89}$, and $35^{118.9}\approx 3.8877\times10^{183}$. These values are so extreme that the limits of floating-point arithmetic are pushed. To achieve better convergence, if one tries further to increase the range of $m$ and $l$, the problem takes the numerical form of an indeterminate $\infty/\infty$. All of this results in a catastrophic loss of precision. 

These numerical limitations explain why series solutions remain impractical for solving problem (\ref{eq:FDE}-\ref{eq:ics}). Our work addresses this issue by developing an approximate closed-form solution for the impulse response of the fractional oscillator that avoids these numerical pitfalls while preserving the fundamental dynamical characteristics of the original system.
\begin{remark}
Note that the impulse response \eqref{eq:imp_resp} can also be written as 
\[
x(t) =  t \sum_{m=0}^\infty \left(-2\zeta (\omega_n t)^{2-\beta}  \right)^m E^{(m+1)}_{2,2+(2-\beta)m} (-\omega_n^2 t^2), 
\]
where
\[
E_{\alpha,\beta}^{(\gamma)}(z) = \frac{1}{\Gamma(\gamma)} \sum_{l=0}^\infty \frac{\Gamma(l+\gamma) z^ l}{l!\Gamma(\alpha l+\beta)}
\]
is the Prabhakar generalisation of the Mittag-Leffler function. In \cite{Garrappa2015}, the authors developed an approach for its implementation in MATLAB (see \cite{garrappa2025mittag}) based on the numerical inversion of its Laplace transform. However, this approach is not directly applicable here as it is restricted for $\gamma\neq 1$ to $\alpha\in(0,1)$.  
\end{remark}
%
%%%%%%
\section{Approximate solution}
\label{subsec:approximate_solution}
%%%%%%%%%%%%%
In \Cref{eq:series_solution}, we have shown that the impulse response can be represented in terms of the double-infinite series (eq. \ref{eq:imp_resp}). We have also demonstrated that when $\beta=0$ and $\beta=1$ (see \Cref{subsec:limiting_case_beta0,subsec:limiting_case_beta1}, respectively), the series converges to the classical impulse solutions consisting of an exponentially decaying term (due to damping $\zeta$) and a sinusoidal function with effective frequency $\omega_d$. The analytical summation of the series was possible as the Gamma functions convert to the factorial functions if $\beta=0$ and $1$. Moreover, we have demonstrated in \Cref{subsec:blow_up_issues} that the series \eqref{eq:imp_resp} suffers numerical instabilities.

\begin{figure}
\centering
\includegraphics[scale=0.45]{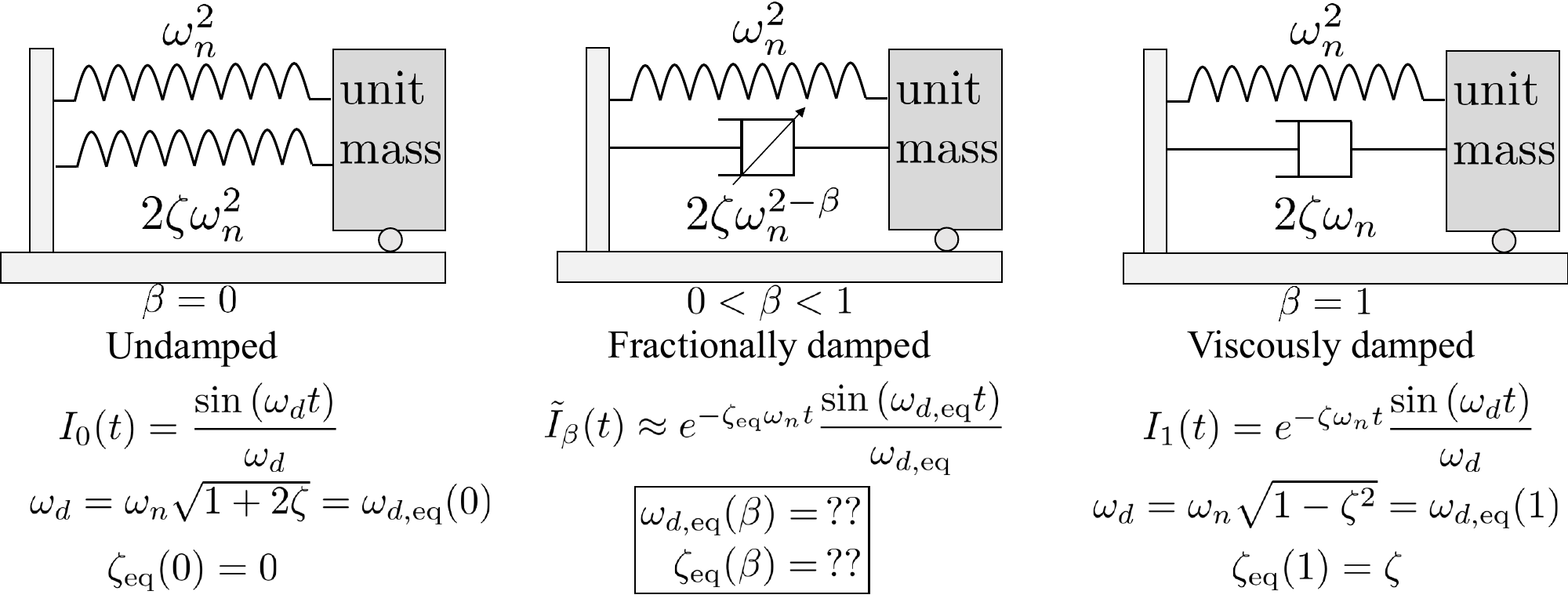}
\caption{Schematic diagram motivating the proposed approximate solution \eqref{eq:approximate_solution}.}
\label{Fig:explain}
\end{figure}

Our main goal is to propose and validate a closed-form approximation of \eqref{eq:imp_resp} that is consistent with the classical cases $\beta=0$ and $1$. From \Cref{Fig:explain}, we can see that $\beta=0$ corresponds to the undamped oscillator with an amplified spring stiffness, whereas $\beta=1$ corresponds to the classical oscillator with viscous damping. Thus, it is intuitive that the oscillator with $\beta$ somewhere in between $0$ and $1$ will exhibit the combined behaviour of a spring and a viscous dashpot. Therefore, the solution for a fractional oscillator can also be approximated to consist of:
\begin{enumerate}[(i), topsep=0pt]
\item an exponentially decaying term with an equivalent damping $\zeta_{\text{eq}}$ as a function of $\beta$ and $\zeta$, and
\item  sinusoidal term with an equivalent damped frequency $\omega_{d,\text{eq}}$ as a function of $\omega_n$, $\beta$ and $\zeta$.
\end{enumerate} 
For these reasons, we propose the following closed-form approximation
\begin{equation} \label{eq:approximate_solution}
\tilde{I}_\beta(t) = \exp\left(-\zeta_{\text{eq}}\omega_n t \right) \frac{\sin \left({\omega}_{d,\text{eq}} t\right)}{{\omega}_{d,\text{eq}}}, 
\end{equation}
where the equivalent damping ratio $\zeta_{\text{eq}}$ and equivalent damped frequency ${\omega}_{d,\text{eq}}(\beta) $ are defined as 
\begin{align}
\zeta_{\text{eq}} (\beta) &:= \zeta \beta^{(0.95 - 0.85 \beta)}, \\
\omega_{d,\text{eq}}(\omega_n,\zeta,\beta) &:= \omega_n \sqrt{1 + 2\zeta - \zeta(2 + \zeta)\beta^{(2.24 - 0.63\beta)}}.
\end{align}
Note that $\tilde{I}_0(t) = I_0(t)$ and $\tilde{I}_1(t) = I_1(t)$, so the approximate impulse response coincides with the solution for $\beta=0$ and $\beta=1$, as required.

In contrast to the impulse response series \eqref{eq:imp_resp}, the proposed approximation \eqref{eq:approximate_solution} consists of the product of an elementary exponential and sine function. Since these functions are implemented in standard numerical libraries using optimised and stable algorithms, the approximation avoids the pitfalls associated with Gamma function computations. The proposed closed-form approximation is more reliable for computations and also faster.  

In the next Subsections (\ref{sec:equivalent_frequency} and \ref{sec:equivalent_damping}), we will motivate the definitions for $\zeta_{\text{eq}}$ and $ {\omega}_{d,\text{eq}}.$ Afterwards, in \Cref{sec:error_analysis}, we will perform an error analysis between the impulse response \eqref{eq:imp_resp} and its approximation \eqref{eq:approximate_solution}.
%%%%%%%%%%%%%
\subsection{Equivalent Frequency of Fractional Oscillators}
\label{sec:equivalent_frequency}
%%%%%%%%%%%%%
In \Cref{subsec:limiting_case_beta0,subsec:limiting_case_beta1}, we defined $\omega_d$ for the limiting cases as 
\begin{equation}
\omega_d = \left\{
\begin{aligned}
  \omega_n \sqrt{1 - \zeta^2} &\text{ for } \beta = 1, \\
         \omega_n \sqrt{1 + 2\zeta}  & \text{ for } \beta = 0.
\end{aligned} \right.
\end{equation}
For $\beta\in [0,1],$ the apparent frequency $\omega_d$ corresponds to the magnitude of the imaginary part of the root of the characteristic equation (see \Cref{eq:laplace_image} in \Cref{appendix:derivation_impulse_response})
\begin{equation}
s^2 + 2\zeta\omega_n^{2-\beta} s^{\beta} + \omega_n^2 = 0. 
\end{equation}
Numerical methods are required for the evaluation of $\omega_d$ as the roots of this transcendental equation cannot be determined analytically for $\beta\in(0,1)$. We used the in-built function `\textit{fsolve}' of MATLAB \cite{matlab_fsolve}, which employs the trust-region dogleg algorithm to solve the equations.

\begin{figure}
\centering
\includegraphics[scale=0.45]{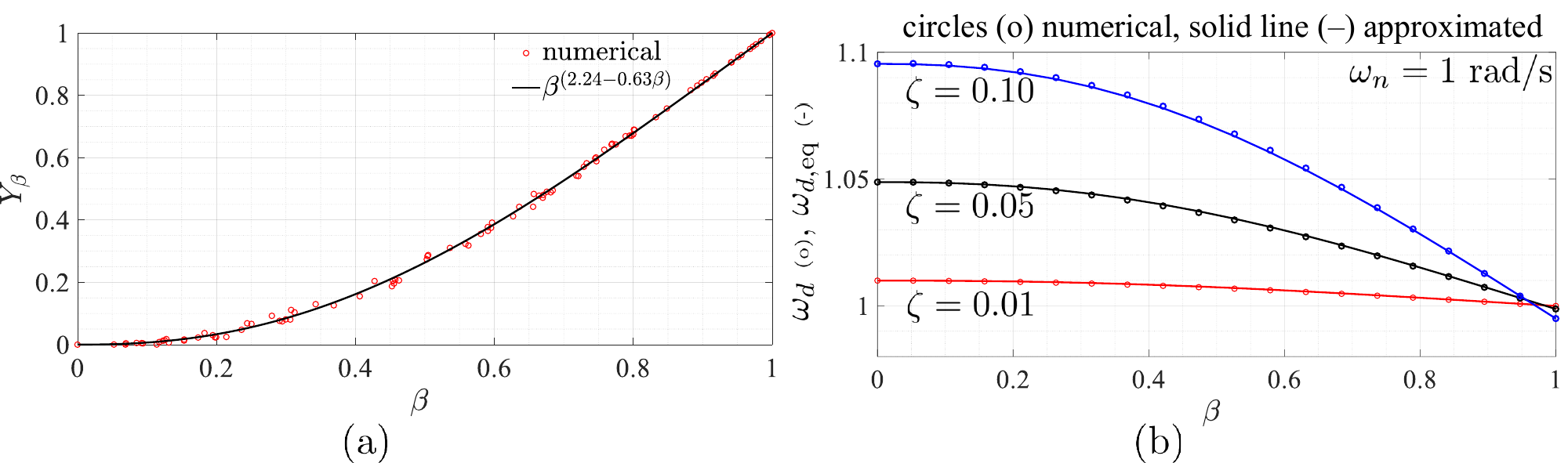}
\caption{(a) Relation between $Y_{\beta}$ defined in \cref{eq:yb} and $\beta$; (b) $\omega_d$ as a function of $\beta$ for different values of $\zeta$.}
\label{Fig:yb}
\end{figure}

It is interesting to note that the fractional damper is an element that has the properties of a viscous damper and a spring depending on the value of $\beta$. For example, when $\beta=1$, it is a purely Newtonian viscous dashpot, and in that case the damped frequency $\omega_d=\omega_n\sqrt{1-\zeta^2}\leq \omega_n.$ When $\beta=0$, there is no damping in the system, and the fractional dashpot behaves like a spring such that the natural frequency of the system becomes $\omega_d=\omega_n\sqrt{1+2\zeta}\geq \omega_n.$ When $0<\beta<1$, there is a partial spring nature in the fractional dashpot, which would increase the system stiffness (and hence there will be an increase in frequency), and the remaining partial viscous damping component would attempt to damp (reduce) the frequency. Let us define this effective frequency as $\omega_{d,\text{eq}}$. From the above explanation, it is obvious that $\omega_n\sqrt{1+2\zeta}\geq\omega_{d,\text{eq}}\geq \omega_n\sqrt{1-\zeta^2}$ and $\omega_{d,\text{eq}}$ decreases gradually from $\omega_n\sqrt{1+2\zeta}$ to $\omega_n\sqrt{1-\zeta^2}$ with increasing $\beta$ from 0 to 1, as can also be seen from the numerical results shown in \Cref{Fig:yb}. Here, we also refer to two interesting papers \cite{naber2010linear, zarraga2019analysis}, where variations of $\omega_d$ with $\beta$ were studied, and similar trends were noted.

To obtain an equivalent expression ${\omega}_{d,\text{eq}}$, we define $M_{\beta}:=(\omega_d/\omega_n)^2$. Note that $M_0 = 1+2\zeta$ and $M_1 = 1-\zeta^2.$ To determine $M_{\beta}$ as a function of $\beta$ and $\zeta$, we introduce a variable $Y_{\beta}$ as follows
\begin{equation}\label{eq:ybeta}
Y_{\beta}=\dfrac{M_0-M_{\beta}}{M_0-M_1}, 
\end{equation}
which scales the transition between $\beta =0$ and $\beta=1.$ We plot the scatter of $Y_{\beta}$ with $\beta$ for randomly generated 10,000 statistically independent and uniformly distributed samples of $\omega_n,\zeta,\beta$. The $\beta$ varies between 0 and 1, $\zeta$ varies from 0.001 to 0.15 (covering the practical ranges), and $\omega_n$ varies from 1 to 10 rad/s (covering the practical range for a variety of civil engineering structures). After conducting several numerical experiments, we observed that the regression model of the form $Y_{\beta}=\beta^{(A_0-A_1\beta)}$ fits well to the scatter, where $A_0$ and $A_1$ are regression parameters. Note that this form captures the power function transition observed in the numerical data. The fractional calculus is based on the generalization of exponential kernels by power-law-type kernels, and, therefore, it appears natural to have a power function here. The regression parameters are obtained using the `\textit{fitnlm}' command in MATLAB \cite{matlab_fitnlm}, which uses the Levenberg-Marquardt least squares approach (considering $A_0=A_1=1$ as initial guesses). From these experiments, we conclude that $A_0=2.2380$ (with a 95\% confidence interval $[2.2320, 2.2441]$) and $A_1=0.6316$ (with a 95\% confidence interval $[0.6217, 0.6415]$). Hence, we have
\begin{equation}\label{eq:yb}
Y_{\beta}=\beta^{(2.24-0.63\beta)}.
\end{equation}
Consequently, using \cref{eq:ybeta,eq:yb}, we get that the apparent frequency $\omega_d$ is approximated by the equivalent frequency $\omega_{d,\text{eq}}$ as follows
\begin{equation}
\omega_{d,\text{eq}}(
\beta)=\omega_n\sqrt{1+2\zeta-\zeta(2+\zeta)\beta^{(2.24-0.63\beta)}}, \quad \beta \in [0,1].
\label{eq:wd}
\end{equation}
\Cref{Fig:yb}(a) shows that the proposed \cref{eq:yb} fits reasonably well with the numerically obtained data for 100 randomly generated test samples. \Cref{Fig:yb}(b) shows a reasonable match between numerically obtained $\omega_d$ and proposed $\omega_{d,\text{eq}}$ for three different damping levels and $\omega_n=1$ rad/s. Note that the expression matches with the limiting cases $\beta = 0$ and $\beta=1.$
%
%%%%%%%%%%%%%
\subsection{Equivalent Damping of Fractional Oscillator}
\label{sec:equivalent_damping}
%%%%%%%%%%%%%
We already know that the damping ratio is $\zeta$ if $\beta=1$ (pure viscous damping) and $0$ if $\beta =0$ (undamped system). However, when $\beta\in(0,1)$, the system is fractionally damped. Therefore, we need to define an equivalent viscous damping ratio $\zeta_{\text{eq}}$, as a function of $\zeta$ and $\beta$, to capture the amplitude decay rate of the fractionally damped oscillator in the classical sense.

\begin{figure}
\centering
\includegraphics[scale=0.45]{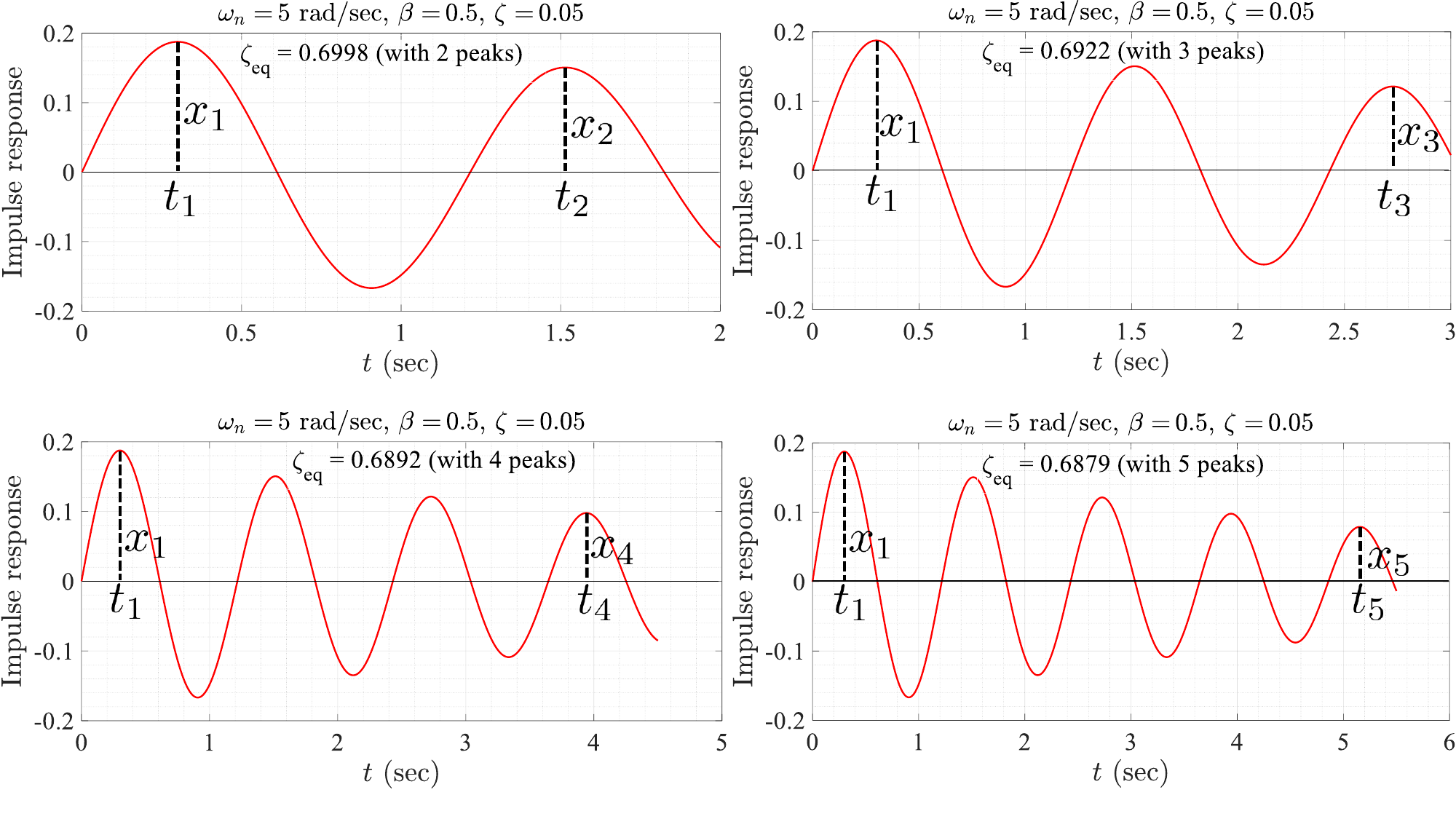}
\caption{Illustration of computation of equivalent damping $\zeta_{\text{eq}}$ via logarithmic decrement method and effect of number of response cycles considered.}
\label{Fig:damping_time}
\end{figure}

\begin{figure}[h]
\centering
\includegraphics[scale=0.35]{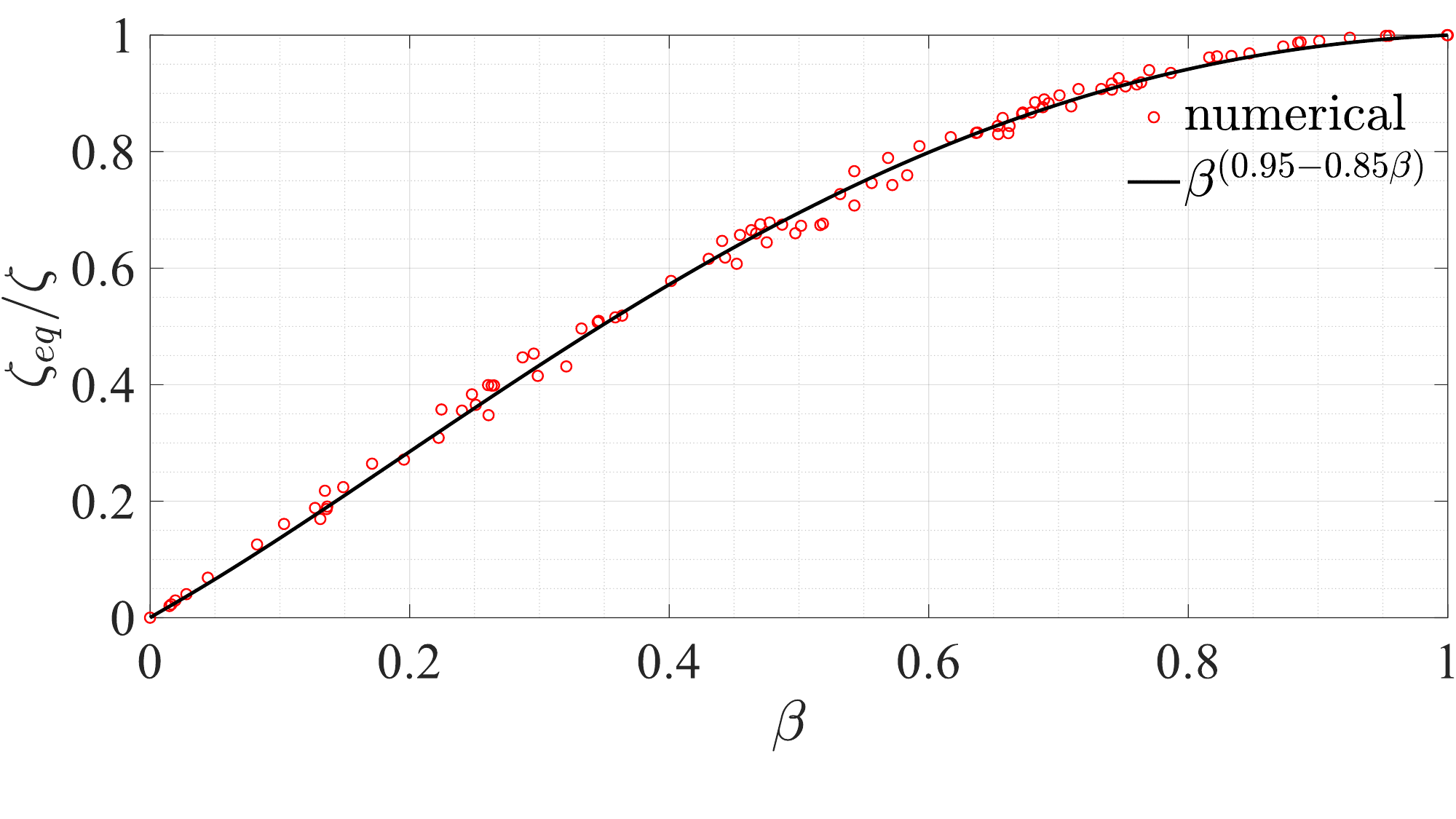}
\caption{Relation between equivalent damping $\zeta_{\text{eq}}$ and $\beta$.}
\label{Fig:damping}
\end{figure}

Using the logarithmic decrement method \cite[Section~2.3]{meirovitch2001fundamentals}, we can estimate the equivalent damping $\zeta_{\text{eq}}$ in the system from a given response signal for $\beta\in (0,1)$. Let the first positive peak in the response $I_\beta(t)$ be $x_1$ and let the $j^{th}$ peak be $x_j$ (refer to \Cref{Fig:damping_time}). If the natural period is $\omega_n$, then, following the logarithmic decrement method, we can write the equivalent damping $\zeta_{\text{eq}}$ as
\begin{equation}
\zeta_{\text{eq}}=\dfrac{\delta}{\sqrt{\delta^2+4\pi^2}},
\end{equation}
where $\delta=\dfrac{1}{j-1}\ln{\left(\dfrac{x_1}{x_j}\right)}$. \Cref{Fig:damping_time} demonstrates the computation of the equivalent damping ratio for an example case with $\beta = 0.5$, $\zeta = 0.05, $ and $\omega_n = 5$ rad/s for $j=2,3,4,$ and 5. It is noted that the rate of response decay reduces slightly with increasing $j$, however, there is no significant difference in the values of $\zeta_{\text{eq}}$ (varies between 0.6879 and 0.6998). Therefore, for the sake of computational efficiency, this study considers only the first two cycles to evaluate $\zeta_{\text{eq}}$.
Now, we will plot $\zeta_{\text{eq}}/\zeta$ as a function of $\beta$ and propose a functional form for $\zeta_{\text{eq}}/\zeta$. For randomly generated 10,000 statistically independent and uniformly distributed samples of $\beta \in (0,1)$, $\zeta \in [0.001, 0.15]$ and $\omega_n \in [1,10]$ rad/s, $\zeta_{\text{eq}}/\zeta$ is calculated as discussed above. We conducted several numerical experiments and fitted the regression model of the form ${\zeta_{\text{eq}}}/{\zeta}=\beta^{(A_0-A_1\beta)}$, where $A_0$ and $A_1$ are regression parameters. As in \Cref{sec:equivalent_frequency}, this form captures the power function observed in the numerical data and satisfies the boundary conditions. From regression analysis, it is noted that $A_0=0.9506$ (with a 95\% confidence interval $[0.9486, 0.9527]$) and $A_1=0.8499$ (with a 95\% confidence interval $[0.8452, 0.8546]$) can be considered as the mean value of the regression parameters. Thus, we have
\begin{equation}
{\zeta_{\text{eq}}}={\zeta}\beta^{(0.95-0.85\beta)}.
\label{eq:damping}
\end{equation}
It is worth mentioning here that the numerical experiments did not indicate any significant relation between $\zeta_{\text{eq}}$ and $\omega_n$. \Cref{Fig:damping} shows that the proposed function \eqref{eq:damping} fits reasonably well with the numerically obtained data for 100 randomly generated test samples.
%
%%%%%%%%%%%%%%%%%
\section{Error Analysis}
\label{sec:error_analysis}
\begin{figure}[h]
\centering
\includegraphics[scale=0.45]{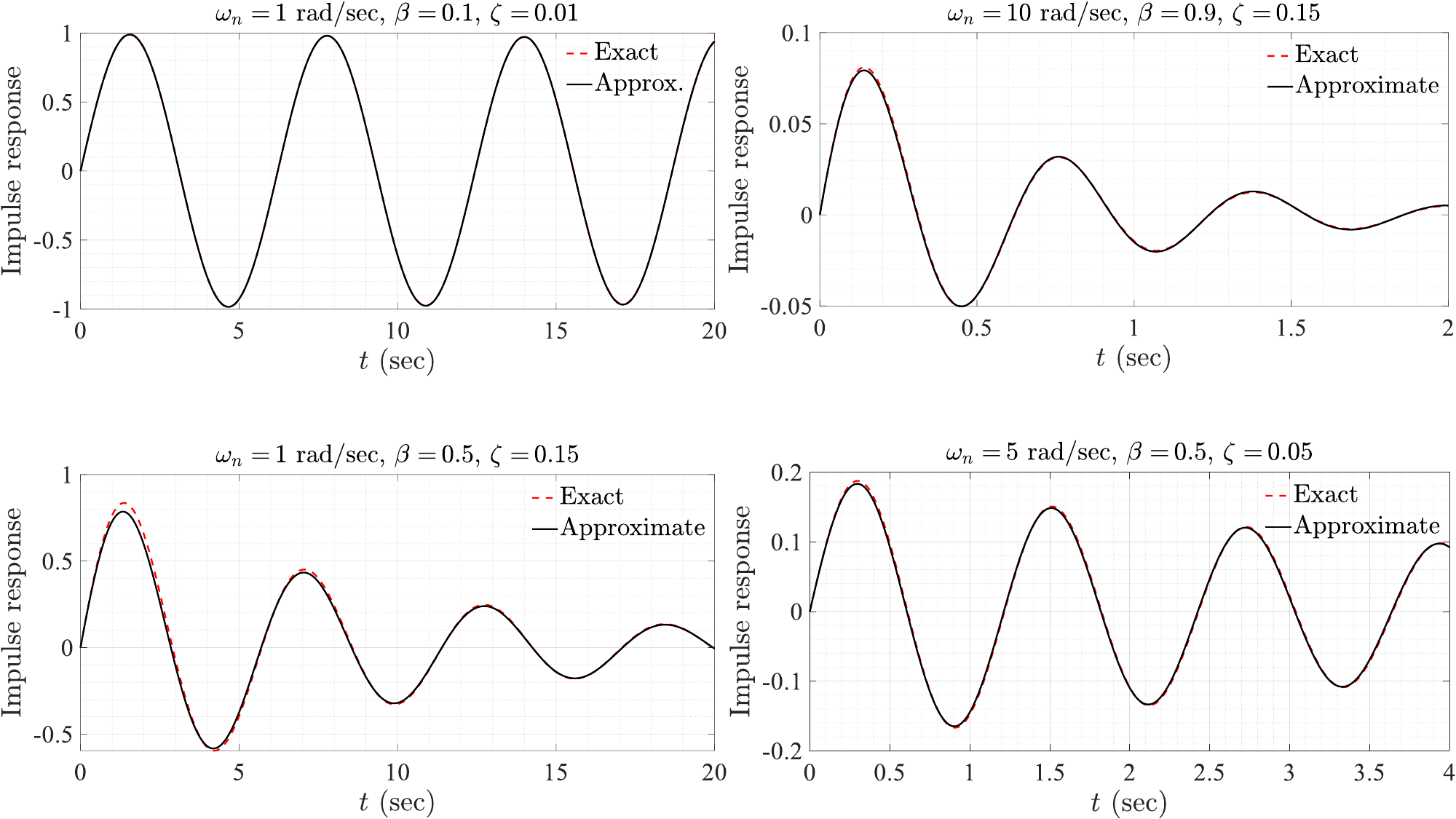}
\caption{Comparison between the series (or exact) solution and the proposed approximate solution for four different sets of $\omega_n, \beta$, and $\zeta$ (values mentioned in the respective titles).}
\label{Fig:error1}
\end{figure}

\begin{figure}[h]
\centering
\includegraphics[scale=0.45]{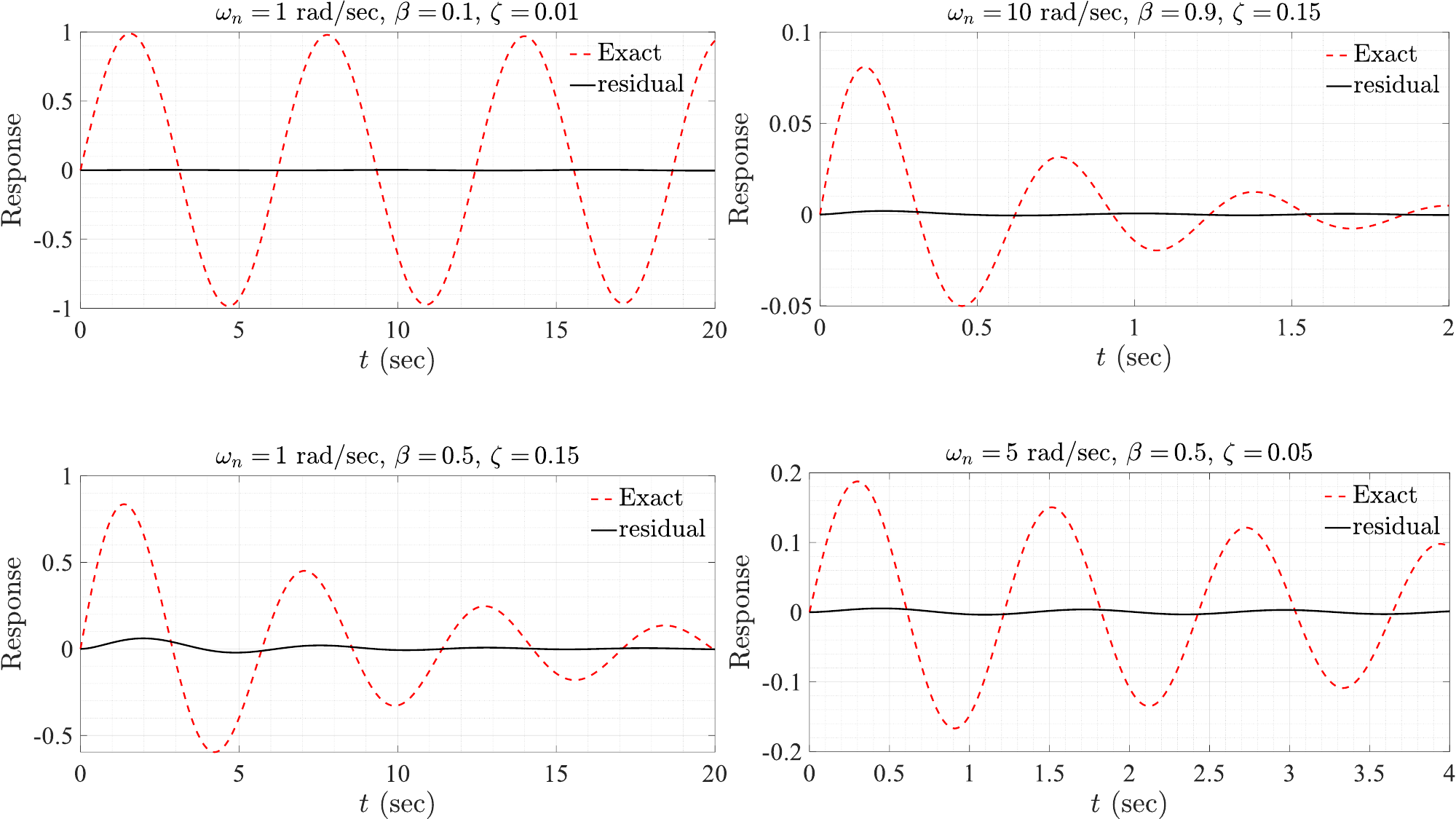}
\caption{Comparison between the series solution and the residuals in the approximate solutions for the four sets of $\omega_n, \beta$, and $\zeta$.}
\label{Fig:error2}
\end{figure}

%%%%%%%%%%%%
In this section, we inspect the performance of the proposed closed-form approximation \eqref{eq:approximate_solution} for the impulse response of the fractional oscillator. The exact series solution \eqref{eq:imp_resp} and the approximate solutions are compared for the following four typical cases:
\begin{enumerate}[(i),topsep=0pt]
\item $\omega_n = 1$ rad/sec, $\beta=0.1$, $\zeta=0.01$, 
\item $\omega_n = 10$ rad/sec, $\beta=0.9$, $\zeta=0.15$,
\item $\omega_n = 1$ rad/sec, $\beta=0.5$, $\zeta=0.15$,
\item $\omega_n = 5$ rad/sec, $\beta=0.5$, $\zeta=0.05$. 
\end{enumerate}
The comparison of responses is shown in \Cref{Fig:error1}, which indicates that the proposed solution is reasonably accurate. Furthermore, we also show the residual in the approximate solution (series solution minus approximate solution) in \Cref{Fig:error2}. The residuals are small enough compared to the exact solution.
%%%%
\subsection{Comparison in Frequency-domain}
\begin{figure}
\centering
\includegraphics[scale=0.45]{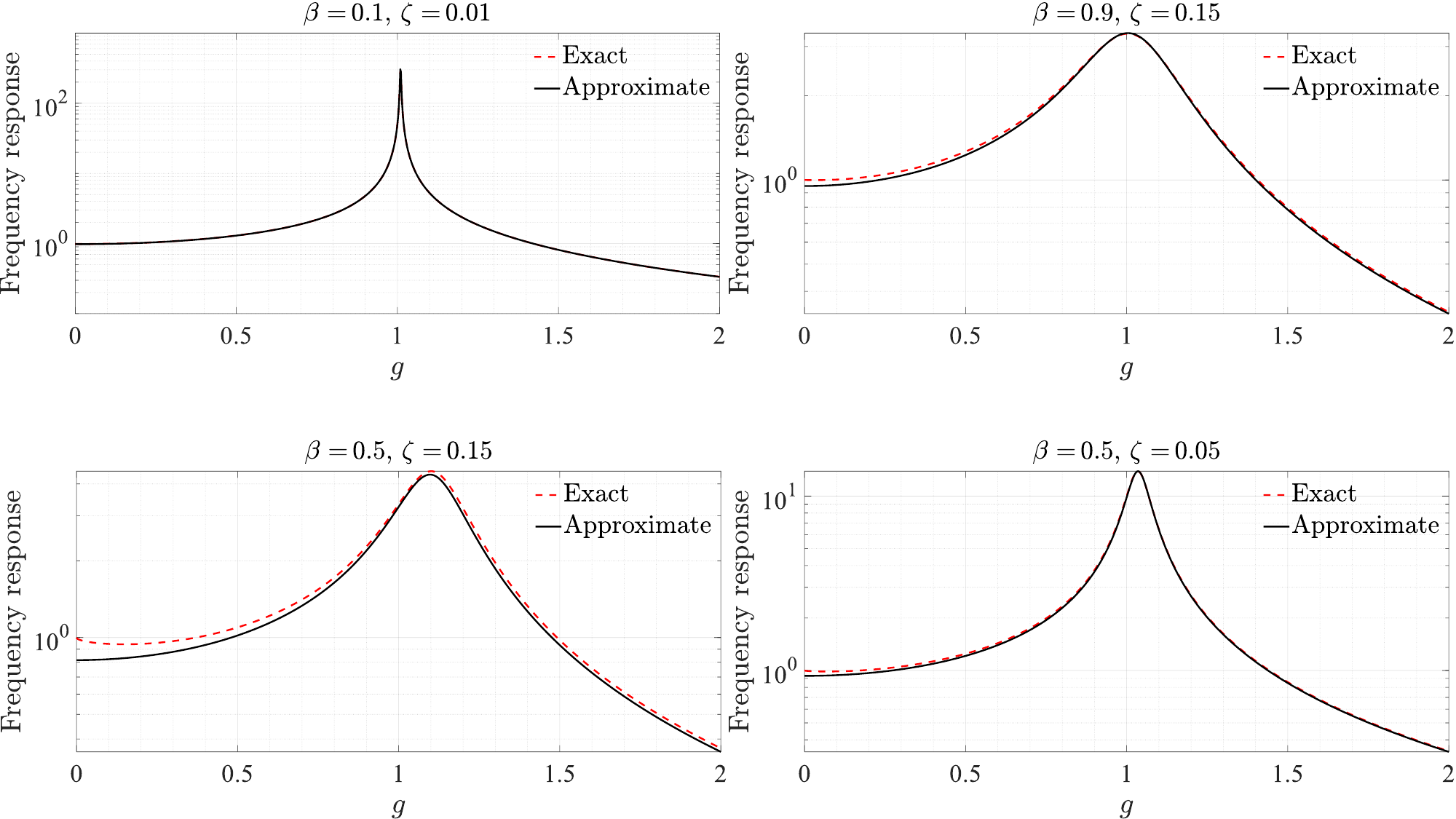}
\caption{Comparison between the frequency domain impulse responses calculated by the exact solution $\lvert h(g)\rvert$ and the approximate solution $\lvert \tilde{h}(g)\rvert$ for four different sets of $\beta$ and $\zeta$.}
\label{Fig:freq_resp}
\end{figure}

%%%5
In addition to the time-domain comparison, we also carry out a comparison of frequency response functions. The frequency response $H(\omega)$ of the exact series solution is obtained by substituting $s=i\omega$ into \cref{eq:laplace_image} (\Cref{appendix:derivation_impulse_response}), where $i=\sqrt{-1}$ and $\omega=$ frequency. This gives
$$
H(\omega) =  \frac{1}{-\omega^2  + 2\zeta\omega_n^{2-\beta} (i\omega)^{\beta} + \omega_n^2},
$$
which can be further simplified in terms of non-dimensional frequency $g=\omega/\omega_n$ as 
\begin{equation}\label{eq:frf}
\omega_n^2H(\omega) := h(g)= \frac{1}{1-g^2  + 2\zeta (ig)^{\beta}}.
\end{equation}
To obtain the frequency response $\tilde{H}(\omega)$ of the proposed approximation \cref{eq:approximate_solution}, we use the definition of the Fourier transform, which leads to
$$
\tilde{H}(\omega) =  {\int_{-\infty}^{\infty}{\tilde{I}_{\beta}(\tau)}e^{-i\omega\tau}\,d\tau} = \frac{1}{\omega_{d,\text{eq}}^2+ (\zeta_{\text{eq}}\omega_n+i\omega)^2},
$$
which in terms of $g$ becomes
\begin{equation}\label{eq:frf_approx}
\omega_n^2\tilde{H}(\omega) := \tilde{h}(g)= \frac{1}{(\omega_{d,\text{eq}}/\omega_n)^2 + (\zeta_{\text{eq}}+ig)^2}.
\end{equation}
The absolute values of the frequency response functions ($\lvert h(g)\rvert$ and $\lvert \tilde{h}(g)\rvert$) are plotted in \Cref{Fig:freq_resp} for the four sets of $\beta$ and $\zeta$ considered previously. A reasonable match is observed between the two solutions. The approximate solution captures the resonance accurately; however, some minor differences are noted in the static response (near $g=0$), especially when the damping $\zeta$ is high and $\beta$ is away from the extremities (0 and 1). This clearly reinforces the robustness of the approximation in the frequency domain as well. 
%%%%%%%%%%
\subsection{Comparison with results of the eigenvector expansion approach}
\begin{figure}[h]
\centering
\includegraphics[scale=0.45]{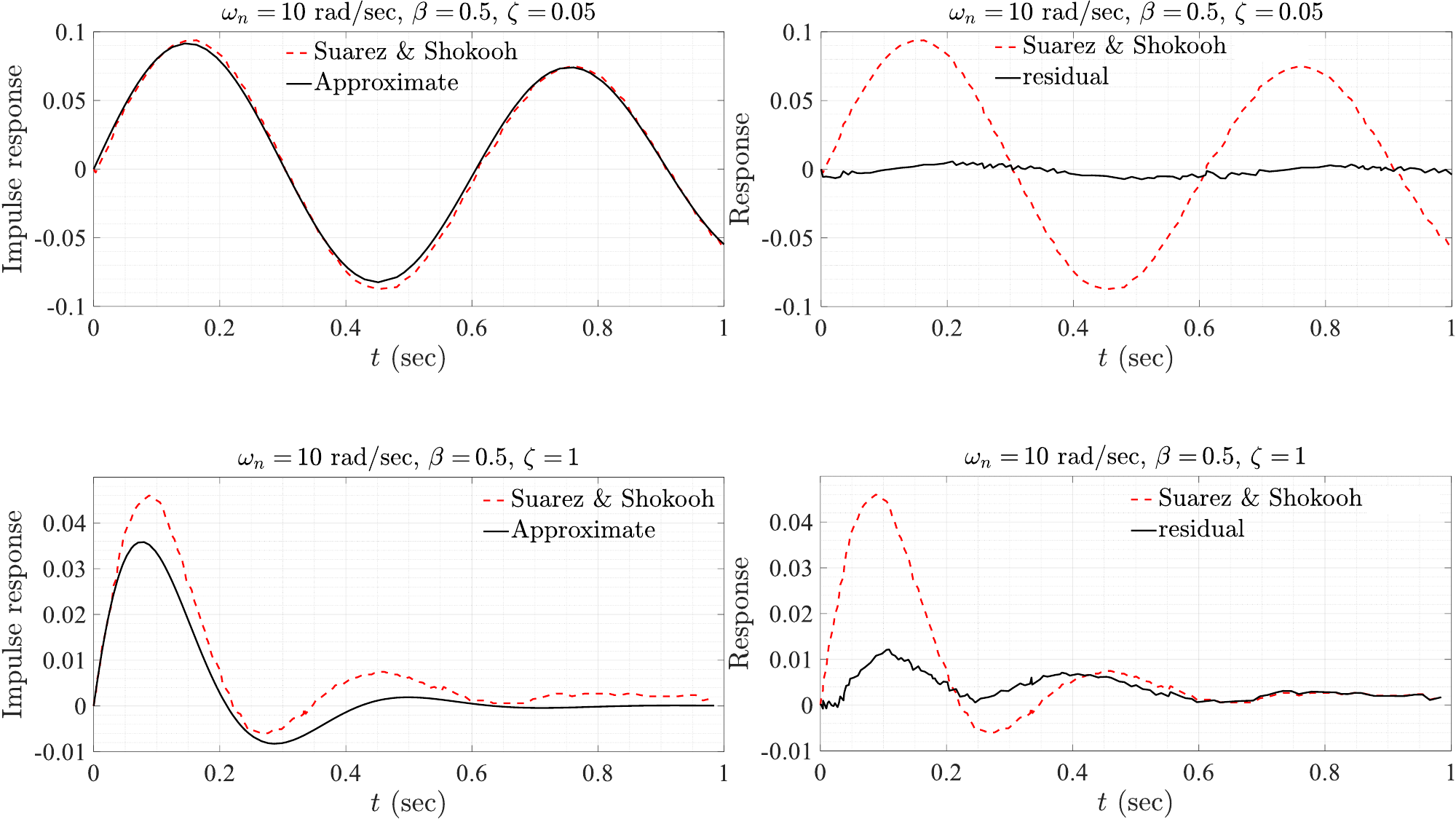}
\caption{The comparison of results reported in Suarez and Shokooh \cite{suarez1997eigenvector} with our approximation.}
\label{Fig:data_comp}
\end{figure}

%%%%%%%%%
Suarez and Shokooh \cite{suarez1997eigenvector} analysed \cref{eq:FDE} in the special case $ \beta=0.5$. They used the eigenvector expansion approach to obtain the impulse response for $\omega_n=10$ rad/s and several damping ratios. We digitised the results corresponding to $\zeta=0.05$ and $1$, contained in their paper, and compared them here with our approximation in \Cref{Fig:data_comp}.
A reasonable match is obtained for $\zeta=0.05$. However, for $\zeta=1$, there is an agreeable qualitative match with some noticeable quantitative deviation because our approximation has been developed for practical values of $\zeta$ which are usually far less than 1. Additionally, the approximated impulse response for $\zeta=1$ depicts an oscillatory nature, which is consistent with the observations of \cite{suarez1997eigenvector, naber2010linear, zarraga2019analysis} where it is discussed that the fractionally damped oscillators do not strictly exhibit the classical overdamped behaviour (i.e., no oscillations when $\zeta=1$).

In a nutshell, all of the above comparisons clearly demonstrate that the proposed approximation is practically acceptable for engineering applications.
%%%%%%%%%%%%%%
\section{Response to Arbitrary Excitation using Convolution}
\label{sec:response}
%%%%%%%%%%%%%%%%
In this section, we demonstrate the application of the developed closed-form approximation $\tilde{I}_\beta$ (see \cref{eq:approximate_solution}) by using it along with the convolution technique (see \cref{eq:sol}) to obtain the responses for any arbitrary excitation.

We consider the same example as in \cite{yuan2017mechanical}, which is given by
$$
\ddot{x}(t)+0.4D^{0.56}x(t)+2x(t)=30\cos{(6t)}.
$$ 
Comparing with the governing equation \eqref{eq:FDE}, we get $2\zeta\omega_n^{2-\beta}=0.4$ (rad/s)$^{1.44}$, $\beta=0.56$, $\omega_n^2=2$ (rad/s)$^2$, and $h(t)=30\cos{(6t)}$ m/s$^2$. Thus, $\omega_n=\sqrt{2}$ rad/s and $\zeta=0.1214$. We can now convolve excitation $h(t)$ with the proposed approximate impulse response \eqref{eq:approximate_solution} and the exact series solution \eqref{eq:imp_resp} to obtain the oscillator response denoted by $x_p(t) = \int_0^t h(t-s) \tilde{I}_\beta(s) \ds$ and $x_s(t) = \int_0^t h(t-s) I_\beta(s) \ds$, respectively.

\begin{figure}[h]
\centering
\includegraphics[scale=0.45]{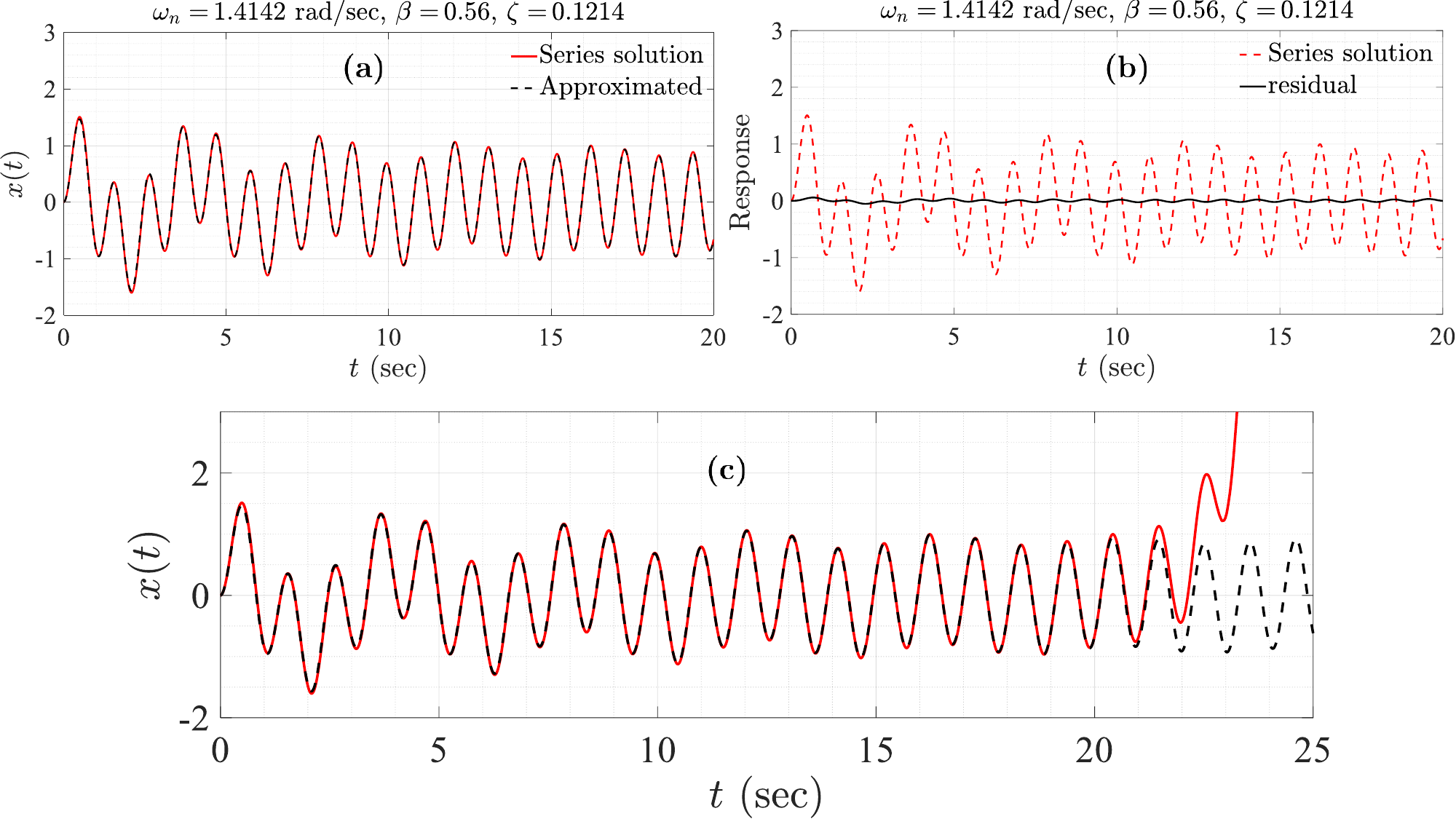}
\caption{Comparison between the responses obtained by convolution on the series solution and on the proposed approximation.}
\label{Fig:compare_resp}
\end{figure}

\begin{figure}[h]
\centering
\includegraphics[scale=0.16]{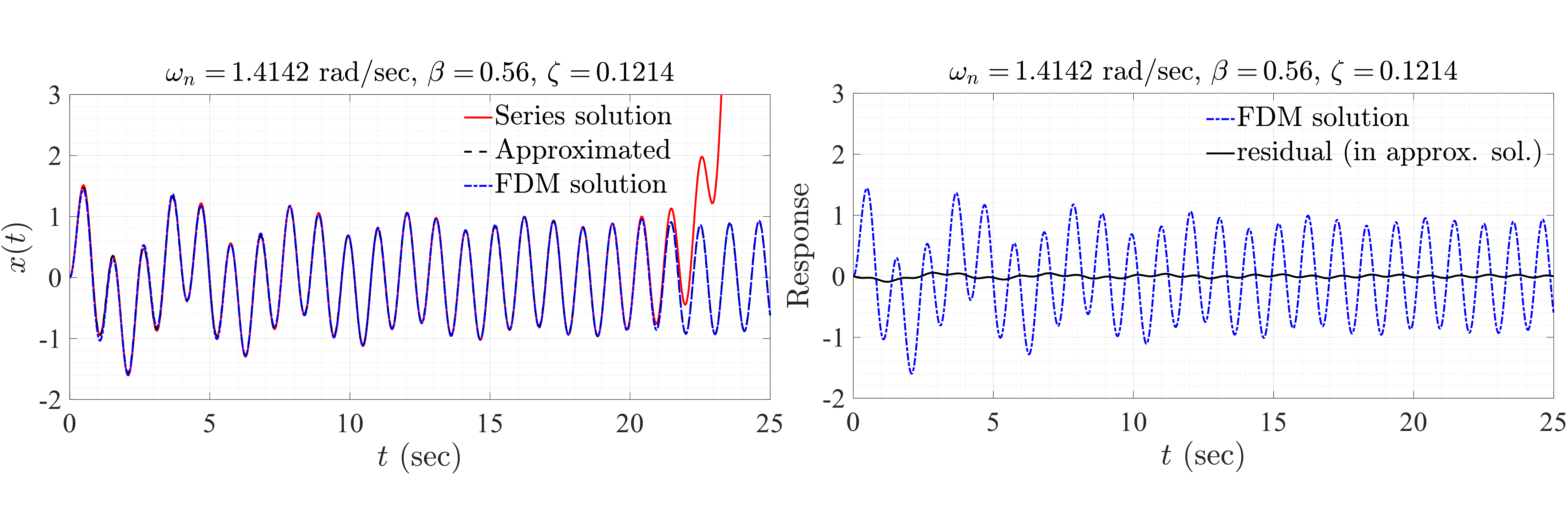}
\caption{Comparison between the responses obtained by convolution on the proposed approximation and the finite-difference solution.}
\label{Fig:compare_resp_fdm}
\end{figure}

The two responses are compared in \Cref{Fig:compare_resp}(a) and the residual $x_s(t)-x_p(t)$ is plotted in \Cref{Fig:compare_resp}(b). A reasonable agreement is obtained between the two responses with negligible residuals throughout. Since the series solution blows up beyond $t=22$s (see \Cref{Fig:compare_resp}(c)), a quantitative comparison could be carried out till $t=22$s only. 
To evaluate the accuracy of the approximated response $x_p(t)$ for $t> 22s$, we compared it with the finite-difference solution denoted by $x_f(t)$. The details of how to obtain the finite-difference (FDM) solution can be found in \Cref{appendix:FDM}. We used a time step $\Delta t=0.005s$ in the implementation of the FDM solution. As can be seen from \Cref{Fig:compare_resp_fdm}, the proposed approximation performs well in the region where the series solution is unable to provide a solution (note the residual $x_f(t)-x_p(t)$).
Thus, we conclude from this section that the proposed impulse approximation leads to reasonably accurate results when convoluted with arbitrary excitation. 
%%%%%%%%%%%%
\section{Conclusion}
\label{sec:conclusion}
%%%%%%%%%%%%
In this work, we considered a fractionally damped oscillator governed by the Caputo derivative of order $\beta\in(0,1)$. The exact solution, expressed as a series of bivariate Mittag-Leffler functions, was shown to recover classical results in the limiting cases of $\beta=0$ (undamped oscillator) and $\beta=1$ (viscous damping). Although the series converges uniformly, its numerical evaluation suffers from instabilities that limit its practical utility. In addition to this, the classical structural dynamics with integer-order derivatives provides a well-known expression for the impulse response of the oscillators; however, no such simple formula exists for the fractional oscillators.

To address the above limitations, we proposed an approximate closed-form equation for the impulse response of fractional oscillators that avoids the computational challenges associated with the Mittag-Leffler function and extends the classical impulse response to the case of fractional oscillators. The equivalent frequency and damping of this approximation were rigorously analysed. It was shown that the proposed form is accurate, and its effectiveness in solving fractionally damped systems under arbitrary external excitation was demonstrated. This confirms its suitability for engineering applications.

Future work could explore the derivation of error bounds between the approximate and exact solutions, as well as the study of the case $\beta\in(1,2)$. 
\appendix
%%%%%%%%
%%%%%%%
\section{Derivation of impulse response}
\label{appendix:derivation_impulse_response}
%%%%%%%%%
We consider the fractional-order differential equation
\begin{equation}\label{eq:prob}
    \ddot{x}(t) + 2\zeta\omega_n^{2-\beta} D^\beta x(t) + \omega_n^2 x(t) = 0,
\end{equation}
with initial conditions $x(0) = 0$ and $\dot{x}(0) = 1$. Taking the Laplace transform,
%and using $ \mathcal{L}[D^\beta x(t)](s) = s^{\beta} \mathcal{L}[x](s) - x(0)=s^{\beta} \mathcal{L}[x](s),$
we obtain
\begin{equation}\label{eq:laplace_image}
    \mathcal{L}[x](s) =  \frac{1}{s^2  + 2\zeta\omega_n^{2-\beta} s^{\beta} + \omega_n^2 }.
\end{equation}
Recall that for $\abs{z}<1$ the geometric series gives $\frac{1}{1-z} = \sum_{k=0}^\infty z^k.$ Put $z = - 2\zeta\omega_n^{2-\beta} s^{\beta-2} -\omega_n^2 s^{-2}$. 
We have $\abs{z} <1$ for sufficiently large $\abs{s}.$
Hence, we obtain 
\[
\mathcal{L}[x](s) = \frac{1}{s^2} \sum_{k=0}^\infty (-1)^k \left(2\zeta\omega_n^{2-\beta} s^{\beta-2} + \omega_n^2 s^{-2}\right)^k. 
\]
Next, employing the binomial theorem $(x+y)^k =  \sum_{m=0}^k \binom{k}{m} x^{m}y^{k-m}$, we expand the $k$-th power and get 
\[
\mathcal{L}[x](s) = \sum_{k=0}^\infty (-1)^k \sum_{m=0}^k  \binom{k}{m}\left(2\zeta\omega_n^{2-\beta} \right)^m \omega_n^{2(k-m)} s^{\beta m -2k-2}. 
\]
As $
\mathcal{L}^{-1} \left\{ s^{-\gamma} \right\} = \frac{t^{\gamma - 1}}{\Gamma(\gamma)}$ for $\Re(\gamma) > 0$, we have 
\begin{equation}\label{eq:sol1}
x(t) = \sum_{k=0}^\infty (-1)^k \sum_{m=0}^k  \binom{k}{m}\left(2\zeta\omega_n^{2-\beta} \right)^m \omega_n^{2(k-m)} \frac{t^{2k + 1 - \beta m}}{ \Gamma(2k + 2 - \beta m) }. 
\end{equation}
We can easily check that, $x(0)=0$ and $\dot{x}(0)=1.$ 

Using $ \sum_{k=0}^\infty \sum_{m=0}^k c_{k,m} = \sum_{m=0}^\infty \sum_{k=m}^\infty c_{k, m} $
and putting $l = k-m$ give 
\[ 
 \sum_{k=0}^\infty \sum_{m=0}^k c_{k,m} =  \sum_{m=0}^\infty \sum_{l=0}^\infty c_{l+m, m}. 
\]
Then, \eqref{eq:sol1} can be written as  
\begin{equation}\label{eq:sol2}
x(t) =  \sum_{m=0}^\infty \sum_{l=0}^\infty (-1)^{l+m}   \binom{l+m}{m}\left(2\zeta\omega_n^{2-\beta} \right)^m \omega_n^{2l} \frac{t^{2 l+2m + 1 - \beta m}}{ \Gamma(2l+2m + 2 - \beta m) }. 
\end{equation}
Using the bivariate Mittag-Leffler function defined in  \cref{eq:mlf}, the solution can be compactly written as 
\begin{equation}\label{eq:sol2b}
x(t) = t  E_{(2,2-\beta),2}(-(\omega_n t)^2,-2\zeta(\omega_n t)^{2-\beta}).     
\end{equation}
%%%%%
\section{Limiting case \texorpdfstring{$\beta \searrow 0$}{}}
\label{appendix:beta_0}
%%%%%
For $\beta =0,$ we have from \eqref{eq:sol1}  that
\[
x(t) = \sum_{k=0}^\infty (-1)^k \frac{t^{2k + 1 }}{ (2k + 1 )! } \sum_{m=0}^k  \binom{k}{m}\left(2\zeta \omega_n^2 \right)^m (\omega_n^2)^{k-m}. 
\]
Using the binomial theorem, we get
\begin{equation}\label{eq:sol_beta_0_eq1}
x(t) = \sum_{k=0}^\infty (-1)^k \left((1 + 2\zeta) \omega_n^2\right)^{k} \frac{t^{2k + 1 }}{ (2k + 1 )! } = \frac{\sin\left(\sqrt{1+2\zeta}\omega_n t\right)}{\sqrt{1+2\zeta}\omega_n},
\end{equation}
since $\sin(x) = \sum_{n=0}^{\infty} \frac{(-1)^n}{(2n+1)!} x^{2n+1}.$
This function $x(t)$ is the solution to the second-order differential equation
\begin{equation} \label{eq:beta_0}
    \ddot{x}(t) + (1+2\zeta) \omega_n^2 x(t) = 0, \quad t>0, 
\end{equation}
with initial conditions $x(0) = 0$ and $\dot{x}(0) = 1$, i.e., it solves problem (\ref{eq:prob}) with $\beta=0$ as $D^0x(t) = x(t)-x(0)=x(t). $ Note that taking the Laplace transform of \eqref{eq:beta_0} gives
\begin{equation} \label{eq:laplace_transform_beta_0}
    \mathcal{L}[x](s)  = \frac{1}{s^2 +  (1+2\zeta) \omega_n^2}. 
\end{equation}
Using $\mathcal{L}^{-1} \left[ \frac{1}{s^2 + a^2} \right] = \frac{1}{a}  \sin(at),$ we get \eqref{eq:sol_beta_0_eq1}, which can be written in terms of the effective natural frequency $\omega_d = \sqrt{1+2\zeta}\omega_n$ as 
\begin{equation}\label{eq:sol_beta_0_eq2}
x(t) = \frac{\sin\left(\omega_d t\right)}{\omega_d} \quad \text{for } \beta =0. 
\end{equation}
Note that $\omega_d$ equals the magnitude of the imaginary part of the root of the denominator of \eqref{eq:laplace_transform_beta_0}. 
%%%%%%
\section{Limiting case \texorpdfstring{$\beta \nearrow 1$}{}}
\label{appendix:beta_1}
%%%%%%
In this section, we will show that taking $\beta =1$ in \eqref{eq:sol2} or \eqref{eq:sol2b} results in the solution
\begin{equation} \label{eq:sol_beta1}
 x(t) = \frac{1}{\omega_d} e^{-\zeta\omega_n t} \sin(\omega_d t), \quad \omega_d = \omega_n \sqrt{1 - \zeta^2}, \text{ for } \beta =1,  
\end{equation}
to the second-order differential equation
\[
\ddot{x}(t) + 2\zeta\omega_n \dot{x}(t) + \omega_n^2 x(t) = 0, \quad t>0,
\]
with initial conditions $x(0) = 0$ and $\dot{x}(0) = 1$, i.e., solving problem (\ref{eq:prob}) with $\beta=1$. First note that 
\begin{equation} \label{eq:laplace_trans_beta_1}
    \mathcal{L}[x](s)  = \frac{1}{s^2 + 2\zeta\omega_n s + \omega_n^2}. 
\end{equation}
As $\zeta < 1$, the roots of the denominator are
$
    s_{1,2} = -\zeta\omega_n \pm i \omega_d. %\quad \text{where} \quad \omega_d = \omega_n \sqrt{1 - \zeta^2}.
$
Rewriting the denominator of \eqref{eq:laplace_trans_beta_1} in terms of these roots results in 
$
    s^2 + 2\zeta\omega_n s + \omega_n^2 
    %= (s + \zeta\omega_n - i\omega_d)(s + \zeta\omega_n + i\omega_d) 
    = (s + \zeta\omega_n )^2 + \omega_d^2.
$
Using the standard inverse Laplace transform identity $
    \mathcal{L}^{-1} \left[ \frac{1}{(s + a)^2 + b^2} \right] = \frac{1}{b} e^{-a t} \sin(bt)$ gives \eqref{eq:sol_beta1}.

Now, we will show that 
\begin{equation} \label{eq:formula_beta1}
    E_{(2,1),2}(-\omega_n^2 t^2,-2\zeta \omega_n t) =  e^{-\zeta \omega_n t} \frac{\sin(\omega_d t)}{\omega_d t}, \quad t > 0. 
\end{equation}
From \cref{eq:mlf}, we have that 
$$
E_{(2,1),2}(-\omega_n^2 t^2,-2\zeta \omega_n t) = \sum_{m=0}^\infty \sum_{l=0}^\infty (-1)^{l+m}   \binom{l+m}{l}\left(2\zeta \right)^m \omega_n^{2l+m} \frac{t^{2l+m}}{ (2l+m + 1)! }. 
$$
Using the Taylor series $\sin(x) = \sum_{n=0}^{\infty} \frac{(-1)^n}{(2n+1)!} x^{2n+1}$ and $e^x = \sum_{n=0}^{\infty} \frac{x^n}{n!}$, and substituting $\omega_d = \omega_n \sqrt{1-\zeta^2}$, we get 
$$
e^{-\zeta \omega_n t} \frac{\sin(\omega_d t)}{\omega_d t}  =  \sum_{m=0}^\infty  \sum_{l=0}^\infty \frac{(-1)^{l+m} \omega_n^{2l+m} (1-\zeta^2)^l  \zeta^m t^{2l+m}}{(2l+1)!m!}. 
$$
Let $ 2l + m + 1  =N $ with $N\in \NN$. This defines a relationship between $ l $ and $ m $ given by $m = N - 2l - 1.$ Since $ m \geq 0 $ and $ l \geq 0 $, the valid values of $ l $ are $0 \leq l \leq \left\lfloor \frac{N-1}{2} \right\rfloor,$ where $\left\lfloor \cdot \right\rfloor$ is the floor function. Hence, we have 
\[
E_{(2,1),2}(-\omega_n^2 t^2,-2\zeta \omega_n t) = \sum_{N=1}^\infty \frac{t^{N-1}}{N!} \sum_{l=0}^{\lfloor \frac{N-1}{2} \rfloor} (-1)^{N - l - 1} \omega_n^{N - 1} (2 \zeta )^{N - 2l - 1}  \binom{N - l - 1}{l}
\]
and 
\[
e^{-\zeta \omega_n t} \frac{\sin(\omega_d t)}{\omega_d t} = \sum_{N=1}^\infty t^{N-1} \sum_{l=0}^{\lfloor \frac{N-1}{2} \rfloor} \frac{(-1)^{N - l - 1} \omega_n^{N - 1} (1 - \zeta^2)^l  \zeta ^{N - 2l - 1}}{(2l + 1)! (N - 2l - 1)!}.
\]
We see that 
\begin{equation} \label{eq:es1} 
E_{(2,1),2}(-\omega_n^2 t^2,-2\zeta \omega_n t) = \sum_{N=1}^\infty t^{N-1} \frac{(-\zeta\omega_n)^{N-1}}{N!}  2^ {N-1}\sum_{l=0}^{\left\lfloor \frac{N-1}{2} \right\rfloor}
  \binom{N-l-1}{l} \frac{ 2^{-2l}} { (-\zeta^2)^l}
\end{equation}
and (using the binomial of Newton)
\begin{align}
e^{-\zeta \omega_n t} \frac{\sin(\omega_d t)}{\omega_d t} & = \sum_{N=1}^\infty t^{N-1} \frac{(-\zeta\omega_n)^{N-1}}{N!} \sum_{l=0}^{\left\lfloor \frac{N-1}{2} \right\rfloor} \left(1-\frac{1}{\zeta^2}\right)^l   \binom{N}{2l+1} \nonumber\\
&= \sum_{N=1}^\infty t^{N-1} \frac{(-\zeta\omega_n)^{N-1}}{N!}  \sum_{l=0}^{\left\lfloor \frac{N-1}{2} \right\rfloor}   \binom{N}{2l+1} \sum_{k=0}^l \binom{l}{k} \left(-\frac{1}{\zeta^2}\right)^k \nonumber\\
&= \sum_{N=1}^\infty t^{N-1} \frac{(-\zeta\omega_n)^{N-1}}{N!}   \sum_{k=0}^{\left\lfloor \frac{N-1}{2} \right\rfloor}   \left(-\frac{1}{\zeta^2}\right)^k \sum_{l=k }^{\left\lfloor \frac{N-1}{2} \right\rfloor} \binom{l}{k} \binom{N}{2l+1}. \label{eq:es2} 
\end{align}
Now, from \cite[Theorem~A]{GLICKSMAN1967} or \cite[Eq.~(3.121)]{Gould1972}, we have
\[
\sum_{l=k }^{\left\lfloor \frac{N-1}{2} \right\rfloor} \binom{l}{k} \binom{N}{2l+1} = 2^{N-1-2k} \binom{N-1-k}{k}. 
\]
Substituting this expression into \cref{eq:es2}, and replacing $k$ by $l$, yields \cref{eq:es1}.
Therefore, we conclude that \eqref{eq:formula_beta1} is valid. 
%%%%%%%%%%
\section{Finite Difference Approach}
\label{appendix:FDM}
For $x\in \Cont^1([0,T])$ satisfying $x(0)=0$, the Gr\"{u}nwald-Letnikov approximation for the Caputo fractional derivative  $D^\beta x(t)$ is given by (see \cite[Section~7.1 and 7.5]{podlubny1998fractional})
\begin{equation}
D^\beta x(t)= \lim_{\Delta t \to 0} \frac{1}{(\Delta t)^\beta} \sum_{j=0}^{\left\lfloor\frac{t}{\Delta t}\right\rfloor} W_{\beta j} \, \, x(t-jh), 
\end{equation}
where 
%\( [\cdot] \) represents the integer part and 
$W_{\beta j}$ are weights defined as 
\[
W_{\beta j} = (-1)^j \binom{\beta}{j}, \quad j = 0, 1, 2, \ldots
\]
To determine the values of $W_{\beta j}$, one can use the following recurrence relation
\[
W_{\beta 0} = 1; \quad W_{\beta j} = \left( 1 - \frac{\beta + 1}{j} \right) W_{\beta,j-1} \quad j = 1, 2, 3, \ldots. 
\]
Hence, we get the following first-order finite-difference form for \cref{eq:FDE}
\[
\dfrac{x(t_{i+1}) -2x(t_{i}) + x(t_{i-1})}{(\Delta t)^2} +\dfrac{2\zeta\omega_n^{2-\beta}}{{(\Delta t)}^{\beta}}\sum_{j=0}^{i} W_{\beta j}\, x(t_{i-j}) + \omega_n^2 x(t_i) = h(t_i).
\]
Rearranging the terms leads to the following recursive relation
\begin{equation}
x(t_i) = \frac{(\Delta t)^2 h(t_i) + 2 x(t_{i}) - x(t_{i-1}) - 2\zeta (\omega_n\Delta t)^{2-\beta} \sum_{j=1}^{i} W_{\beta j}\, x(t_{i-j})}{1 + 2\zeta (\omega_n\Delta t)^{2-\beta} + (\omega_n \Delta t)^2}, \quad i = 2, 3, 4, \ldots
    \label{Fdm_recursive_relation}
\end{equation}
where
\[
  x(0) = 0, \quad \dot{x}(0) \approx \frac{x(t_1)-x(0)}{\Delta t} =0,
\quad \text{for sufficiently small}\, \Delta t.
\]
\begin{acknowledgement}
The suggestions of Professor Arnaud Deraemaeker (Universit\'e Libre de Bruxelles, Belgium) on frequency domain comparisons are well acknowledged.
\end{acknowledgement}

\bibliography{references}
\bibliographystyle{unsrt} %
\end{document}